\input amstex
\documentstyle{amsppt}

\def\rightrhead{FACTOR = QUOTIENT, UNCOUNTABLE BOOLEAN ALGEBRAS  \hfil}

\document
\headline{\tenrm\rightrhead}
\def\restriction{|}

\topmatter
 \title FACTOR  = QUOTIENT, UNCOUNTABLE BOOLEAN ALGEBRAS, \\
NUMBER OF ENDOMORPHISM AND WIDTH \endtitle
\author Saharon Shelah \endauthor
\thanks{Partially supported by Basic Research Fund of the Israeli Academy of Sciences. Publ. No. 397}\endthanks

\address{The Hebrew University of Jerusalem, Israel,
Rutgers University, New Brunswick, NJ USA, MSRI, Berkeley, Calif, U.S.A.}\endaddress
\date{June 22, 1991}\enddate
\abstract{
We prove that assuming suitable cardinal arithmetic, if $B$ is a Boolean algebra every homomorphic image of which is isomorphic to a factor, then $B$ has locally small density. We also prove that for an (infinite) Boolean algebra $B$, the number of subalgebras is not smaller than the number of endomorphisms, and other related inequalities. Lastly we deal with the obtainment of the supremum of the cardinalities of sets of pairwise incomparable elements of a Boolean algebra.}
\endabstract
 \endtopmatter

\def\mpr{\medskip\par}
\def\mpi{\medskip\par{\it Proof.\; \;}}

We show in the first section:

\mpr{\bf 0.1 Conclusion}. \; {\sl It is consistent, that for every
 $Q=F$ Boolean algebra $B^*$, for some $n<\omega,
\{x:B^*\restriction x \text{ has a density} \leq \aleph_n
\}$ is dense (so $B^*$ has no independent subset of power
$\aleph_n)$.
}

Where:

\mpr{\bf 0.2 Definition}. \;{\sl A \  $B$.$A.$\  is \ $Q=F$ 
(quotient equal factor) if: every homomorphic image of $B$
is isomorphic to some factor of $B$ i.e $B\restriction a
\,\text{for some}\,  a\in B.$
}
\medskip
 The ``consistent'' is really a derivation of the
conclusion from a mild  hyphothesis on cardinal arithmetic
(1.2). The background of this paper is a problem of Bonnet
whether  every $Q=F$ Boolean algebra is superatomic.

Noting that : \ ``$B \restriction x$ has density $\leq
\aleph_n$'' is a weakening of  ``$x$ is as atom of
$B'$'', we see that {\bf 0.1} is relevant. 

The existence of non trivial example is proved in R. Bonnet, S. 
 Shelah [2].

M. Bekkali, R. Bonnet and M. Rubin [1] characterized all
interval Boolean algebras with this property.

 In the second section we give a more abstract version.
 In a paper in  preparation , Bekkali and the author use
theorem {\bf 2.1} \ to show that every $Q=F$ tree Boolean
algebras are superatomic.

In the third section we deal with the number of
endomorphism (e.g. $aut(B)^{\aleph_0}\le end(B)$) and in the fourth
with the width of a 
Boolean algebra.

We thank D.~Monk for detecting an inaccuracy in a previous version. 

\medskip
 {\bf Notation}
 
$B$ \ \ \ \ \ \ denote a Boolean algebra.

$B^+$\ \ \ \ \  is the set of non zero members of B. 

$B\restriction x$ \ \ \ \ where $x\in B^+$ is
$B\restriction \{y:y\leq x\}$.

$\text{comp}\,B$ \ is the completion of $B$, so it is an
extension of $B$.  

Remember : if $B_1$ is a subalgebra of
$B_2$, $B$ is complete, $h$\,a homormorphism from $B_1$ to
$B$, then $h$ can be extended to a homomorphism from
$B_2$ to $B$. 

$id_A$ is the identity function on $A$.

${}^\gamma\alpha=\{\eta:\eta$ is a sequence
of lenghts $\gamma$ of ordinals $<\alpha\}$.

${}^{\gamma>}\alpha=\bigcup_{\beta<\alpha}{}^\beta\alpha$.
  
 \mpr{ \S \bf 1 Maybe every ``quotient equal factor" $BA$
  has locally small density }. 

 \mpr{\bf 1.1 Hypothesis}. \; {\sl
 \roster
 \item $B^*$ is a 
$Q=F$ Boolean algebra.
\item  for each $n<\omega$, $B^*$ has a factor $B_n$
s.t.:
$0<x\in B_n \Rightarrow \text{density}
(B_n \restriction x)\geq \aleph_n$.   
\endroster 
}

\mpr{\bf 1.2\  Hypothesis}. \; {\sl For 
$\alpha\geq \omega$ we have
$2^{|\alpha|}>\aleph_\alpha$. 
}

\mpr{1.3\  Desired Conclusion}. \; {\sl Contradiction.

We shall use {\bf 1.1} all the times, but {\bf 1.2} only in
{\bf 1.17}.
}

\mpr{\bf 1.4 \  Definition}. \; {\sl $K^*_\lambda$ is the class of
Boolean algebra such that:
$$(\forall x\in
B^+)[\text{density}(B \restriction x)=\lambda].$$
}

\mpr{\bf 1.5\  Claim}. \; {\sl If $B$ is atomless, $x\in B^+$ then
for some $y$, $0<y\leq x$,and infinite cardinal $\lambda$
we have $B \restriction y \in K^*_{\text{density} (B
\restriction x)}$. 
}

\mpr{\bf 1.6\  Claim}. \; {\sl
If $B\in K^*_\lambda$, \quad $B\subseteq B'\subseteq
\text{comp} (B)$ then,
 $$B'\in K^*_\lambda.$$
}

\mpr{\bf 1.7\  Claim}. \; {\sl If $B\in
K^*_\lambda$,$\aleph_0\leq\mu<\lambda$, $\mu$ regular 
{\it then} some subalgebra $B'$ of $B$ is in
$K^*_\mu$.  If in addition $B$ is a $Q=F$ algebra, {\it
then} some homomorphic image $B''$ of $B$ is in
$K^*_\mu$.
}
 
\mpi Choose by
induction on $i<\mu$, $B_i\subseteq B$, $| B_i | \leq \mu$,
$[i<j\Rightarrow B_i\subseteq B_j]$ \ such that:

if $x\in B^+_i$ then there is
$y(x,i)\in B_{(i+1)}$ satisfying:
\roster
 \item $0<y(x,i) \leq x$,
\item for no $z\in B_i$, \  $0<z\leq y(x,i)$.
\endroster
(possible as for each $i$ and $x \in B^+_i\,
density(B_i \restriction x) = \lambda)$. Let
now  $B'=\cup_{i < \mu} B_i$; it is a subalgebra
of $B$.  Now $x\in (B')^+\Rightarrow density
(B' \restriction x)=\mu$ as on the one hand
$|B'|\leq \sum_{i<\mu} |B_i| \leq \mu \times
\mu=\mu$ implies $density(B' \restriction x) \leq
\mu$ for every $x \in B'$ and on the other hand if $x \in
(B')^+, A \subseteq B' \restriction x, |A| < \mu$ then for
some $i < \mu, A \subseteq B_i $, hence
$y(x,i) \in (B')^+$ wittness $A$ is not dense in $B'
\restriction x$,
Now, if $B$ is a $Q=F$ Boolean Algebra, then
$id_{B'}$ can be extended to a  homomorphism $h'$ from $B$
into $comp(B')$, so  $h'(B)$ is as required.

\mpr{\bf 1.8\  Conclusion}. \; {\sl $\{\lambda: \lambda$  regular
and $B^*$ has a factor in 
 $K^*_\lambda\}$ is an initial segment of
$\{\aleph_\alpha : \aleph_\alpha \,\text{regular}\}$.
}

\mpr{\bf 1.9\  Definition}. \; {\sl 
\roster 
\item $\alpha(*)$ is 
minimal such that for no $\lambda\
>\aleph_{\alpha(*)}$ does $B^*$ has a factor in
$K^*_\lambda$. Let $\kappa^*=\kappa(*)=:|\alpha(*)|$.
\item Let for $\alpha <\alpha^*$, $b_\alpha\in B^*$
be such that $B^* \restriction b_\alpha \in
K^*_{\aleph_{\alpha +1}}$. 
\item Let $J_\alpha$ be the ideal\ $\{b\in B^*:  B^*
\restriction b\in K^*_{\aleph_{\alpha+1}}\}$.
 \endroster
}

\mpr{\bf 1.10\  Definition}. \; {\sl $B[^{\omega >}\lambda]$ is 
the Boolean algebra generated freely by $\{x_\eta:\eta \in
{}^{\omega >}\lambda\}$ except $x_\eta \leq x_{\eta
\restriction m}$ \ $(m\leq \, lg \,(\eta), \eta
\in{}^{\omega >}\lambda)$. }

\mpr{\bf 1.11\  Claim}. \; {\sl 
\roster
\item For $\alpha < \beta< \alpha^*$, $J_\alpha \cap
J_\beta=\{0\}$ and $J_\alpha$ is an ideal.
\item For no $B'  \subseteq B^*$ and proper ideal $I$
of $B'$, $[x\in B' \backslash I \Rightarrow density(B'
\restriction  x/I) >\aleph_{\alpha(*)}]$. \endroster
}

\mpi

(1) Trivial.

(2) By {\bf 1.5} for some $\beta > \alpha(*)$ and proper
ideal $J$ of $B'$ (of the form $\{x \in B':x-b \in I\})$ we
have $B'/I\in K^*_{\aleph_\beta}$ so there is a
homomorphism $h$ from $B^*$ into  $comp(B'/I)$ extending
$x\mapsto x/I\, (x\in B')$.  So $B^*$ has a factor
isomorphic to $Rang\,h$, but this Boolean algebra is in
$K^*_{\aleph_\beta}, \aleph_\beta>\aleph_{\alpha(*)}$.  So
by {\bf 1.7} we get contradiction to the choice of
$\alpha(*)$.

\mpr{\bf 1.12\  Definition}. \; {\sl 1) $I^*=:\{x\in B^*:  \,
\bigcup_{(\alpha<\alpha(*))} J_\alpha$ is dense  below
$x$\}.

2) For $A \subseteq \alpha(*): I^*_A=:\{x\in B^*:\bigcup
\,\{J_\alpha : \alpha\in A\}$ is dense below $x\}$.
}

 \mpr{\bf 1.13\  Claim}.
{\sl  \roster
\item $I^* $, $I^*_A$ are ideals of B.
 \item $A\subseteq
B\Rightarrow I^*_A \subseteq I^*_B \subseteq I^*$.
\endroster }  

\mpr{\bf 1.14\  Claim}. \; {\sl For every $A\subseteq
\alpha(*)$, there are $c_A$,
$h_A$ such that:  \roster 
\item $c_A \in B^*$ and $I^*_A$ is dense below $c_A$,
\item $h_A$ is a homomorphism from $B^*$ onto $B^*|
c_A$, 
\item $h_A \restriction I ^*_A$ is one to one,
\item If $B \restriction x\, \cap
I^*_A=\{0\}$ then $h_A (x)=0$,  
 \item $B^* \restriction c_A$ is
a subalgebra of  the completion of the subalgebra
 $\{h_A (x) : x \in I^*_A\}$. 
\endroster
}

\mpi
Let $ba(A)$ be $$I^*_A \cup \{1-x:x\in I^*_A\},$$
this is a subalgebra of $B^*$.

Let $h_1$ be a homomorphism from $B^*$ to
$\text{comp}\bigl(b \, a (A)\bigr)$ extending $id_{ba (A)}$
and $h_1(x)=0$ \  if\  $B^*\restriction x\cap
I^*_A=\{0\}$.

Now $h_1 (B^*)$ is a quotient of $B^*$ hence there is an
isomorphism $h_2$ from $h_1 (B*)$ onto some
$B^*\restriction c_A$.

Let $$h_A=h_2\circ h_1$$
so (1), (2), (3), (4), (5) holds.

\mpr{\bf 1.15 \ Claim}. \; {\sl Let $A\subseteq \alpha(*)$
\roster 
\item If $x \in \bigcup_{\alpha \in A} J_\alpha$ then
$h_A(x) \in \bigcup_{\alpha \in A} J_\alpha$.
\item  If $x \in I^*_A \setminus \bigcup_{\alpha
\in A} J_\alpha$ then $h_A(x) \not\in \bigcup_{\alpha \in
A} J_\alpha$.
 \item $\bigcup_{\alpha \in A} h_A  (J_\alpha)$
is dense and downwared closed in $B^* \restriction c_A$.
\item $c_A \in I^*_A$.
\item If $x \in J_\alpha$ and $\alpha \in A$ then $h_A(x)
\in J_\alpha$. \endroster
}

\mpi \roster
\item By (5) of {\bf 1.14}.
\item is easy too.
\item is easy too.
\item is easy too.
\item is easy too.
\endroster

\mpr{\bf 1.16\  Claim}. \; {\sl We can find $x_\eta \in B^*$ for
$\eta \in{} ^{\omega >}\lambda$ where $\lambda
=2^{\kappa(*)}$
such that:
\roster
\item $m<lg(\eta) \Rightarrow 0<x_\eta
\leq x_{(\eta \restriction m)}$.
\item If $\bigwedge^k_{e=1}\nu_e \ntrianglelefteq
\eta$, $\eta \in{} ^{\omega >}\lambda$, $\bigwedge_e \nu_e
\in {}^{\omega }\lambda $ then $$B^* \vDash x_\eta
\not\leq \bigcup^k_{e=1} x_{\nu_e}.$$
\endroster
}

\mpi Now we can choose $\langle A^0_\eta :\eta
\in{}^{\omega >}\lambda \rangle$ which is a family of
subsets of $\kappa(*)$ such that any non trivial Boolean
combination of then has cardinality $\kappa(*)$.

Let for $\eta \in {}^{\omega >}\lambda$ \quad
$A_\eta=^{def} \bigcap_{e\leq lg(\eta)} A_{\eta\restriction
e}$.  Let $$x_{<>}=c_{A_{<>}}=h_{A_{<>}} \,
(1_{B^*}). \tag a$$ $$x_{<i>}=h_{A_{<>}}
(c_{A_{<i>}})=h_{A_{<>}}h_{A_{<i>}}\  (1_{B^*})
\tag b$$
\noindent
and generally,

$$x_{<{i_0, i_1}, \dots, i_{n-1}>} = h_{A_{<>}}
h_{A_{<i_0>}} \, h_{A_{<i_0,i_1>}} \dots
h_{A_{<{i_0, i_1},\dots ,i_{n-1}>}} \ (1_{B^*}).$$ 

We prove {\bf (a)} by induction of $lg\eta$.

\mpr
The reader may check

\mpr{\bf 1.17\  Final Contradiction}. \; {\sl $\{x_\eta :\eta \in
{}^{\omega >}\lambda\}$ from {\bf Claim 1.16} contradict by
{\bf 1.11(2)} and the choice of $\aleph_{\alpha(*)}$,
because
$\lambda=2^{\kappa(*)}=2^{|\alpha(*)|}>\aleph_{\alpha(*)}$.

[of course $\aleph_{\alpha(*)} \leq | B^*|$]

Actually, we have prove more.
}

\mpr
{\bf 1.18  Remark}. \;  (1) So we have in {\bf 1.17} prove that
if set theory is as in {\bf Hypothesis 1.2}, then there is
no Boolean algebra as in {\bf 1.1}, hence proving {\bf
Conclusion 0.1}

(2) Note: if  {\bf 1.2}, any $Q=F$ 
Boolean algebra has no factor $\Cal P(\omega)$.

\mpr{\bf  \S 2 $Q=F$ Boolean algebras:\  a general
theorem }.

\mpr{\bf 2.1 \ Theorem}. \; {\sl Suppose:
 \roster
\item $B^*$ is a $(Q=F)$ Boolean algebra.
\item $N$ is a family of (non zero) members of $B^*$ (the
``nice'' elements).
\item $\kappa$ a cardinal $(\geq \aleph_0) ,\ \  
\langle K_\alpha : \alpha <\kappa \rangle$ a sequence.
\item $K_\alpha$ is family of Boolean algebras closed under
isomorphism and for $\alpha \not=\beta$ we have $K_\alpha
\cap K_\beta=\emptyset$.
\item for every $\alpha$ some  factor of $B^*$ is in
$K_\alpha$.
\item if $x\in (B^*)^+, \, \,  (B^*|x)\in K_\alpha$
 then for some $y\leq x,\, B^* \restriction
y \in K_\alpha,\, y\in N$. 
\item if $x_1, x_2 \in N,
\quad B^*\restriction x_1 \in  K_{\alpha{_1}}, \quad
B^*\restriction x_2\in K_{\alpha{_2}}, \quad
\alpha_1\not=\alpha_2$ then $x_1\cap x_2=0$.
\item if $x\in B^*, N'\subseteq N$
{\it then}

$(\alpha)$ \ \  for some $y_1,\cdots,y_n \in N'$,for every
$z\in N'$ we have $B^*\vDash x\cap z \subseteq
y_1\cup\cdots\cup y_n$

or

$(\beta)$ \qquad for some $y\in N': y\leq x$.

Actually we use (8) only for $N'$ of the form $\{y \in N
:(\exists \alpha)[\alpha \in A\, \&\, B^* \restriction y\in
K_\alpha\}$. 
\item if $x<y\in B^*, 
B^*\restriction x\in K_\alpha$ and $B^*\restriction y \in
K_\beta$,  then\  $\alpha=\beta$. 
\endroster
{\it Then} the Boolean algebra $ B \left[ ^{\omega
>}(2^\kappa)\right]$  can be embedded into $B^*$,
remember $^{\omega >}(2^\lambda)$ is 

the tree $\{\eta : \eta $ a finite sequence of
ordinals $ < 2^\kappa\}$ and Def 1.10.
}

\mpr{\it Proof of theorem \  2.1}.

Let for \  $\alpha < \kappa$ ,\, $Y_\alpha=:\{y\in N:B^*|
y\in K_\alpha\}$.

For $A \subseteq \kappa$ let us define

$I_A=$\,the ideal generated by\ $Y_A= \bigcup_{\alpha\in A}
Y_\alpha=\{y\in N : B^*|y\in K_\alpha$ for
some $\alpha \in A\}$ and

$J_A=:\{z\in B_\alpha : \text{for every}\,  y\in
Y_A\,\text{we have}\, z\cap y=0\}$.

Clearly  $J_A$ is an ideal.

Now for each $A\subseteq \kappa$ $B^*/J_A$ is a
quotient of $B^*$.  Hence by condition (1) there are $y^*_A
\in B^* $ and an isomorphism $h_A:B^*/J_A\rightarrow
B^*|y^*_A$ onto. 

Let $g_A:B^*\rightarrow B^*/J_A$ be canonical, 
 so $h_A\circ g_A(1_{B^*})=y^*_A$. Let $f_A=h_A\circ g_A$.

Define for $y\in B^*$ the following: $\text{cont} (y)=:
\{\alpha: (\exists y'\leq y)[B^*|y'\in K_\alpha]\}$.

(i.e. the content of $y$). We next prove

$(*)_1$ \quad cont $(y^*_A) \supseteq A$.

\mpi By conditon (5) for each $\alpha \in A$,
there is $x_\alpha \in B^*$, such that:
$B^*\restriction x_\alpha \in K_\alpha$.
 By condition (6) wlog $x_\alpha \in N$,
hence $x_\alpha \in I_A, \text{hence}\,
g_A \restriction(B^*\restriction x_\alpha)$, is one to one,
hence $B^*\restriction x_\alpha \cong B^*\restriction
f_A (x_\alpha)$, hence by condition (4) $B^*\restriction
f_A (x_\alpha) \in K_\alpha ; \text{now}\,
Rangf_A=B^*\restriction y_A^*$, so $\alpha \in \text{cont}
(y^*_A)$. So we have prove $(*)_1$. 

\mpr
$(*)_2$\   cont$(y_A^*) \subseteq A$

\mpi  Suppose $\alpha \in cont(y_A^*), \text{so
there is}\, z\leq y_A^*,\text{ such that}\,
B^*\restriction z\in K_\alpha$. As $f_A$ is a
homomorphism from $B^*$ onto $B^*\restriction y^*_A$, there
is $x\in B^*$ such that  $f_A(x)=z$.

Now the kernel of $f_A$ is $J_A$, and
$B^*\restriction 0\not\in K_\alpha$ so $x\not\in J_A$;
and clearly ($B^* \restriction x)/ J_A$ is in
$K_\alpha$.

Hence by $(*)_3$ below $\alpha\in A$. 

\mpr

$(*)_3$ \  if $ x\in B^* \setminus J_A$ and\ 
$(B^*/J_A)\restriction f_A(x)\in K_\alpha \,\text{then}\,
\alpha \in A$

\mpi We apply condition (8) to $x$ and
$N'=Y_A$.  So one of the following two cases
ocurrs:

\mpr{\bf Case $\alpha$: } \; There are
$n<\omega, \, y_1,\cdots,y_n\in N' $  such that:

$(\forall z \in N') \, x\cap z \subseteq
y_1\cup\ldots\cup y_n$. 

So $ x-(y_1\cup\cdots\cup y_n)\in J_A
\quad \text{(by definition of}\, J_A)$.

Let $x_1=x\cap (y_1\cup\cdots\cup y_n) \, $ hence
$B^*\restriction f_A(x)=(B^*\restriction x)/J_A\cong
(B^*\restriction x_1)/J_A \cong B^*\restriction x_1$,(last
isomorphism as $\wedge_e y_e \in Y_A$ hence $y_1\cup
\cdots\cup y_n \in I_A$ hence $f_A\restriction
(B^*\restriction x_1)$ is one to one). So $B^*\restriction
x_1\in K_\alpha, \, \text{hence by condition (6) for
some}\, x_2\leq x_1$,we have $ x_2\in N \& B^*\restriction
x_2 \in K_\alpha$.

 Let  $B^*\restriction y_e\in K_{\alpha{_e}}
\text{where}\, \alpha_e \in A \text{(by definition  of} \; 
Y_A)$. Clearly $x_2\leq x_1 \leq y_1 \cup\dots\cup y_n$,
so for some $e, y_e\cap x_2\not=0$ hence
$\alpha=\alpha_e\in A$ (by condition (7)) so we get the
trivial desired conclusion.

\mpr{\bf Case $\beta$: } \;  There is $y\in
N', \, y\leq x$.

As $y \in N'= Y_A = \cup_{\beta\in A} Y_\beta$ for some
 $\beta\in A$ we have
$y \in K_\beta$, also $y\in N'\subseteq N$ so
 $J_A \cap (B^*\restriction y)=\{0\}$ so  $(B^*/
J_A) \restriction g_A(y) \in K_\beta$.  Remembering
$(B^* / J_A)\restriction g_A(x) \in K_\alpha$,
as $B^*/J_A\cong B^*\restriction y_A^*$ we get by condition
(9) that $\alpha=\beta$.  So $(*)_3$ hence $(*)_2$ is
proved.

\mpr
Next we prove
 
$(*)_4$ \quad if $B \subseteq A \subseteq \kappa$, and
 cont$ (y)=B$, {\it then}  $cont[f_A\,(y)]=B$.

\mpi

\mpr{\bf inclusion \  $\supseteq$ }

First let $\alpha \in B$, then for some $x \leq
y$,$B^* \restriction x\in K_\alpha$, and by condition (6) $wlog \, x \in N$,  hence\  $x\in I_A\,
(\text{as} \,  \alpha \in B\subseteq A)$ hence
$f_A\restriction (B^*\restriction x)$ is one to one and
onto $B^*\restriction f_A(x)$ so $f_A(x)\leq f_A(y), \,
B^*\restriction f_A(x) \in K_\alpha$,  so $\alpha \in
\text{cont}(f_A(y))$. 

\mpr{\bf inclusion\  $\subseteq$ }

Second let us assume $\alpha \in cont[f_A(y)]$. So (as
$f_A$ is onto $B^*\restriction y^*_A$, and if $f_A(x) \leq
f_A(y)$ then $f_A(x\cap y)=f_A(y), \,x\cap y \leq y)$
there is $x\leq y $ such that $B^*\restriction f_A(x) \in
K_\alpha$.  Now apply condition (8) to $x$ and $Y_A$. So
case $(\alpha)$ or case $(\beta)$  below  holds.

\mpr{\bf Case\  $\alpha$: } There are $n <
\omega, \, y_1,\cdots,y_n \in Y_A$ such that for every \,
$z\in Y_A$ we have \ $x \cap z \subseteq y_1 \cup
\cdots\cup y_n$. 

Hence $x - (y_1 \cup \cdots \cup y_n) \in J_A$ let $x_1
=x\cap (y_1 \cup...\cup y_n)$ so $f_A(x)=f_A(x_1) $ so
$B^*\restriction f_A(x_1)\in K_\alpha $ and of course $x_1
\leq y_1 \cup\cdots\cup y_n, \,  (\text{and}\,x_1 \leq y)$
so $ f_A \restriction(B^* \restriction x_1) $ is one to
one. 

Now $f_A$ is one to one  on $B^* \restriction x_1$ hence
$B^* \restriction x_1 \cong B^*
\restriction  f_A(x_1) \in K_\alpha$.  Now $x_1 \leq x
\leq y$, so $x_1$ wittness $\alpha \in cont(y)$, which is
$B$. 

\mpr{\bf Case  $\beta$: } There is $t\leq x,
\, t\in Y_A$.

Now $t\leq x \leq y, B^*\restriction t \in
K_\beta$ \,for some $\beta \in A$ as $t \in Y_A$.  Now
$f_A$ is one to one on $B^*\restriction  f_A (t) \in
K_\beta$.

Also $f_A(t)\leq f_A(x)$ hence by assumption (9) we have
$\beta=\alpha$.  Also $t\leq x \leq y$ so $t$ wittness
$\beta \in \text{cont}(y)$, so\  $\alpha = \beta\in
\text{cont}(y)=B$ as required.

So we have proved $(*)_4$

\mpr{\it end of proof of theorem 2.1}: \; Let
$\lambda=2^\kappa$.

Let $\langle \Cal U_\eta:\eta \in{}^{\omega >} \lambda)$
be a family of subsets of $\kappa$, any finite Boolean
combination of them has power $\kappa$ (or just $\not=
\emptyset$).

Let $\Cal U^*_\eta=\cap_{e\leq lg\eta}
 \Cal U_{\eta\restriction e}$.  Now define 
 for every 
$\eta \in {} ^{\omega >} \kappa$ and $e\leq
lg(\eta)$ an element $y^e_\eta$ of $B^*$:

$y^e_\eta= ^{def}\,  f_{{\Cal U_{\eta|e}^*}}$ 
$f^*_{\Cal U_{\eta | (e+1)}} \cdots$  $f_{\Cal
U^*_{\eta\restriction (n-1)}}$  $f_{\Cal
U^*_\eta}(1_{B^*})$ and $y^\otimes_\eta=^{def}
y^\circ_\eta	$

Now: (a) prove for each $\eta \in{}^n\kappa$ by downward
induction on $e\in\{0,1,...,n\}$ that 
$cont(y^e_\eta)=\Cal U^*_\eta$ ; for $e=n$ this is
$(*)_1  +(*)_2$ as $y^n_\eta=y^*_{\Cal U_\eta}$;

for $e < n$ (assuming for $e+1$) this is by $(*)_4$.

Next note: (b) if $\nu=\eta{}\,\hat{} <\alpha>$ {\it then}
$y^\otimes_\nu \leq y^\otimes _\eta$

[prove by domnward induction for $ e\in
\{0,1,...,lg\eta\}$ we have $: y^{e+1}_v \leq y^e_\eta$;
remember $f_{\Cal U}$ is order preserving].

Lastly note (c) if $\eta\in {}^{\omega >}\lambda,
n<\omega$, and $\nu_e \in{}^{\omega >}\lambda$ is not
initial segment of $\eta$ for $e=1,\dots,n$ then
$y^\circ_\eta-\cup^n_{e=1} y^\circ_{\nu_e}\not=$; this
follows by (a) and the definition of
$cont(y^\circ_\eta)$.

Now by (a), (b), (c) there is an embedding $g$ from the
subalgebra of $B^*$ which $\{y^\circ_\eta:\eta \in
 {}^{\omega >}\lambda\}$ generates mapping $y^\circ_\eta$ to
$x_\eta$ .


\def\mue{\mu}
\def\<{\langle}
\def\>{\rangle}

\def\restrictedto{\restriction}
\def\restricted to{\restrictedto}
\def\composition{\circ}
\def\union{\cup}
\def\intersection{\cap}

 \mpr{\bf  \S 3 The Number of Subalgebras }.

\mpr{\bf 3.1 Definition}. \;  For a BA  A 
\roster
\item  $Sub(A)$ is the set of subalgebras of $A$.
\item  $Id(A)$ is the set of ideals of $A$.
\item$End(A)$ is the set of endomorphisms of $A$.
\item $Pend(A)$ is
the set of partial endomorphisms (
i.e.  homomorphisms from a subalgebra of $A$\ into $A)$.
\item $Psub(A)$  is the family of subsets of A 
closed under union ,intersection and substruction 
but 1 may be not in it though 0 is [so not neccessarily
closed under complementation].
\item  We  let $sub(A)$, $id(A)$, $end(A)$, $aut(A)$, 
$pend(A)$, $psub(A)$
 be the cardinality of  $Sub(A)$, $Id(A)$, $End(A)$,
 $Aut(A), Pend(A)$ and $Psub(A)$ respectively
\endroster

\par\quad           
In D. Monk [4] list of open  problems appear:

\mpr{\bf PROBLEM 63}. \;   Is there a BA A such that $aut(A) > sub(A)$ ?
             
See [4] page 125  for backgraound.
             
\mpr{\bf 3.2 Theorem}. \; {\sl For a $BA\, A$ we have: 
$aut(A)$\ is not bigger than $sub(A)$.
}
             
\mpr{\bf 3.3 Conclusion}. \; {\sl 
\roster 
\item For a $BA \,A$ we have
$end (A)$ is not bigger than $sub(A)$.
\item For a $BA \, A$ we have \ 
$pend(A)$ is exactly  $sub(A)$.
\item  For a $BA \, A$ and $a$ in $A,$ $0<a<1$ we have
$sub(A) =Max \{ sub(A \restrictedto a), 
sub(A \restrictedto -a) \}$.
\endroster
}

We shall prove it in 3.5
\mpr{\bf Remark}. \;  Of course - $A$ is infinite- we many times
forget to say so.
                 
\mpr{\it 3.4 Proof of the theorem}. \; 
 Let $\mue$ be  $sub(A)$. 

{\bf 3.4A Observation}: $P sub(A)$  has cardinality
$Sub(A)$ [why?  for the less trivial inequality,
$\leq$,  for every $X$ in $P sub(A)$ which is not a subset of
$\{0,1\}$ choose a member $a\not=0,1$ in $X$ and
let $Y_a[X]$ be the subalgebra generated by $\{ x: x \in X,
x  \le   a \}$ let $Z_a[X]$ be the subalgebra generated by
$\{ x: x \in X , x \intersection a =0 \}$;
now $X$ can be reconstructed from 
$\< a,Y,Z \> $ as $\{x\in A:x\cup a\in Y$ and $x-a\in Z\}$. So $|P
sub(A)| \leq |A| \times |Sub(A)|^2+4 
=|Sub(A)|+\aleph_0=|Sub(A)|$. (remember that $|A|\ge\aleph_0$, and
that $|A|\le |Sub(A)|$, as the number of non-atoms of $A$ is $\le
|\{Z_a[A] : a\in A\}|$).      
             
For any automorphism  $f$ of $A$ we shall choose a
finite sequence of members of $Psub(A) $
( in particular of ideals of $A$ ), and this mapping is
 one to one,
thus  we shall finish.
             
Let $J =^{def}  \{ x : x \in A, \hbox{and for every} \;
y\in A \text{ below} \;  x , \hbox{ we have  }f(y)=y \}$;
let $I^*$ be the ideal of $A$ generated by $I$, the set of elements $x$
for which $f(x) \intersection x  =0$.

{\bf Observe} : $x  \in I$ implies $f(x) \in I$ [why?\,
     as $y=f(x), x \in I $ implies $f(x) \cap x=0$
hence $f(y) \intersection y=
        f(f(x)) \intersection f(x) =
        f(  f(x) \intersection x )=f(0)=0$].
             
{\bf Observe}: $J \union I$ is a dense subset of $A$ [if
$x\in B^+$ and there is no $y\in J, 0<y\leq x$ {\it then}
$ wlog \, x\not= f(x)$.  If $x\not\leq f(x)$ {\it then} $z=:
x-f(x)$ satisfies $0<z\leq x, f(z)\cap z\leq
f(x)\cap z=f(x) \cap (x-f(x))=0$; so $z \in I$.  If $x\leq
f(x)$ {\it then} $x<f(x)$; as $f$ is an automorphism of
$B^*$, for some $z$ in $B^*$ we have $f(z)=x$, so  $x<f(x)$
means $f(z) <f(f(z))$ hence $z <f(z)$, and $z'=x-z=f(z)-z$
is in $I$, is $>0$ but $< x$, as required.]
             
Next let $\{ x_i:i<\alpha \} $  be a maximal sequence of
distinct members  of $I$ satisfying : for any $i,j < 
\alpha $ we have  $x_i \intersection  f(x_j) =0$,
let $I_1$ be the ideal generated by $\{ x_i:i < \alpha
\}$ and  let $I_2$ be the ideal generated by $\{f(x_i) :i <
\alpha \} $.

 Clearly  $I_1  \intersection  I_2 = \{ 0 \}$
let $I_0 = \{ y : y \in I^*, \hbox{ and for every }x \in I_1 ,
                y \intersection x =0
                \hbox{ and }f(y) \in  I_1 \} $
and let $I_3$ be the ideal of $A$ generated by $\{ f(x) : x
\in I_2\} $
             
{\bf Observe} : each $I_t$  ($t =0,1,2,3$) is an ideal of
$A$, contained in 
        $I^*$  [why? 
for $t=1$ by choice each $x_i $
           is in $I$, 
$t=0$ by it's
      definition , 
$I_2, I_3$ by their definition and an observation above, i.e., $x \in I \Rightarrow f(x) \in I$]
             
{\bf Observe}: for $t=0,1,2 $ we have: $I_t \intersection
I_{t+1} =\{0\}$
      [why ? for $t=1$, by the choice of
 $\{ x_i:i<\alpha \}$,
      for $t=0$ by definition of $I_0$, lastly
      for $t=2$ applies its definition + $f$ being
an automorphism.]
                         
 {\bf Observe}: $I_0 \union I_1 \union I_2 $ is a dense
subset of 
        $I^*$    
      [why? assume $x$ in $I^*$ but below it there is no non zero member
      of this union , so we can replace it by any non zero 
      element below it; as $x \in I^*$, there is below it 
       a non zero element $y$ with  $y \intersection f(y) =0$
       so wlog $x \intersection f(x) = 0  $;
       why have we not choose $x_\alpha =x$?
     there are two case:

\mpr{\bf {Case 1}}: For some $j<\alpha$,  $x
\intersection f(x_j) $ is not zero

   so there is a non zero element below $x$ in $I_2$.

\mpr{\bf Case 2}: for some $j<\alpha, f(x)
\intersection x_j$ is not zero

    so there is a non zero  $y$ in $I_1 $ below $f(x) $
     hence (as $f$ is an automorphism)
        there is $x' $ below $x$ such that $f(x') =y$
     so wolg $f(x)$ is in $I_1$, but then by its
 definition, either below $x$ there is a member of $I_1$ or $x$ is in $I_0$
so we have finished proving the observation.]

{\bf Observe}: for $t=0,1,2$, $x\in I_t \Rightarrow f(x)\in I_{t+1}$
[check].
                       
Now we define $C_t=C^f_t$ , a member of $Psub(A)$  for
$t=0,1,2$ as follows: $C_t $ is the set $\{x \union f(x) 
: x \in I_t \}$. The closure under the relevant operations
follows  as $I_t$ is closed under them and $f$ is an
isomorphism and  for $x, y$ in $I_t$,  $x \intersection
f(y) = 0$, this is needed for substraction.
             
 Also for every automorphism $f$ of $A$
we assign the sequence 
$\< J,I,I_0,I_1,I_2,I_3, C_0,C_1$, $C_2 \>$
(some reddundancy).
Suppose for $f_1, f_2 \in Aut(A)$ we get the same tuple;
it is enough to show that
their restriction to J and to I are equal
-as the union is dense.
Concerning $J$ this is trivial - they are the identity
on it        
so we discuss $I$, by  an observation above it  is
enough to choeck it for $I_t , t=0,1,2$
but for each $t$, from $C_t$ and $I_t,I_{t+1}$ we can
 [or see below] reconstruct $f \restriction I\sb t$.
             
So we have finished the proof.
             
\mpr{\bf 3.4B Remark}. \; We can phrased this argument as a claim:
 let $I,J$, be ideals of $A$ with intersection $\{0\}$;
for every $f$, a one to one homomorphism from $I$ to $J$
let $X_f$ be the set $\{x \union  f(x) : x \in I \}$
{\it then} ``$f $ mapped to $X_f$"
is a one to one mapping from $HOM(I,J)$ to $Psub(A)$
(the former include $\{ g \restrictedto I : g \in Aut(A),
g$ maps $I$ onto $J \} $ , for which we
use  this). For subalgebra relative to $\cup, \cap$ only $f$ needs not be one to one.
            
\mpr{\it 3.5 Proof of the conclusion from the theorem}. \;
\roster
\item For a given Boolean Algebra $B$ assume $\mue =:
sub(A)$ is $< end(A) $. For any endomorphism $f$ of $A$
we attach  the pair $(Kernel(f), Range (f))$ .
The number of possible such pairs is 
at most 
$id(A)\times  sub(A)$ , 
which is  at most $\mue$ (we are dealing with infinite
$BAs$ and $id(A) \leq psub(A) \leq sub(A)$) so as we assume
$\mue < end(A)$, there are distinct $f_i$ , endomorphisms
of $A$ for $i<\mue^+$ , an ideal $I$ and a subalgebra $R$
of $B$ such that for every $i$ we have
$Kernel (f_i) = I$ and $Range(f_i) =R$.
             
We now define a homomorphism $g_i $ from $B/I$ to $R$ by :
$g_i (x/I) = f_i(x)$  . Easily the definition does not
depend on the representative,
so $g_i$ is as  required and it is one to one
and onto. So $\{ g_i \composition (g_0)^{-1}
: i < \mue^+ \}$ is a set of $\mue^+$ distinct
 automorphisms
of $R$.

So 

$$\mue < aut(R) \leqno (*)_1 $$
\noindent
but, by the theorem 

$$ aut(R) \le  sub(R) \leqno (*)_2$$
\noindent
obviously  
  
$$ sub(R) \le   sub(A)  \leqno (*)_3$$
\noindent
remember 
     
$$ \mue  = sub(A) \leqno (*)_4$$
\noindent
together a contradiction

\item  As  $R \mapsto$  ( identity map on $R$) is 
one to one from $sub(A)$ into $Pend(A)$, obviously $sub(A)
\leq pend(A)$, so we are left with proving the other
inequality.
 Same proof, only the domain is a subalgebra too
and it has an ideal.
 So for every such partial endomorphism $h$ of   $A$ 
we attach   two subalgebras $D_h =Domain(h)$ and $R_h =Range(h)$
an ideal $I_h$ of $D\sb h\  \{ x : x \in D\sb h$ and 
$f(x)=0 \}$ . They are all in $Psub(A)$, so their number is
at most $sub(A)$ and if we fixed them
the amount of freedom we have left is :
 an automorphism of $R$ (and $aut(R) \le sub(R) \le sub(A)$).
\item  Let $B $ be a subalgebra of $A$.  We shall
attached to it several ideals and subalgebras of
$B\restriction a, B\restriction(-a)$ such that $B$ can be
reconstructed from them; this clearly suffice.  Let $C$  
be the subalgebra  $\< B,a\> , C_o=\{x\in C : x\leq a \},
C_1=\{x\in C: x\leq 1-a\}$. The number of possible $C$ is
clearly  the number of pairs $\langle C_o, C_1 \rangle$
which is clearly $sub(A \restrictedto a) + sub(A
\restrictedto -a)$ ; fix $C$. Let $I=: \{ x: x
\in B, x \le  a \}$, it is an ideal of  $C \restrictedto
a$ ; 
 so the  number of such $I$ is at most

$id(C \restrictedto a)\leq pend(C\restrictedto A) \le 
sub(C \restrictedto a) \leq sub(A \restrictedto a)$.

So we can fix it .
Similarly we can fix $J= \{ y: y \in B, y  \le   (1-a) \}$,
now $I$ and $J$ are subsets of $B$ , 
now check : the amount of freedom we have left 
is an isomorphism $g$   from $(C \restrictedto a)/I$
onto $(C \restrictedto -a)/J$
such that $B= \{$ the subalgebra of $C$ generated by 
$I \union J \union \{x \union y :  x \in (C \restrictedto a),
 y \in (C\restricted to  -a)$  and  $g(x/I)=(y/J)\}$.
 \endroster
            
 So  we can finish easily.

\mpr
We originally want to prove $Aut(A)^{\aleph_0}\leq$ sub
$A$ and even 

$|\{f\in End(A):\,f$ is onto $\}|^{\aleph_0}\leq sub A$. 
But we get more:
 intermediate invariants with reasonable
connections.

\mpr{\bf 3.6 Definition}:  For a Boolean algebra $A$
\roster
\item"(0)"  A partial function $f$ from $A$ to $A$ is
everywhere onto if:  

$x\in Dom(f) \& y \in Rang(f) \& y
\leq f(x) \Rightarrow (\exists z) [z \leq x \& f(z)=y]$.
\item"(1)" $End_0(A)= End(A)$.

 $End_1(A)=\{f\in End(A):\, Rang(f)$ include a
dense ideal$\}$.

 $End_2(A)$ is the set of endomorphism $f$ of $A$
which are onto.

$End_3(A)=\{f:$ for some dense ideal $I, f$ is an
onto endomorphism of $I\}$.

$End_4(A)=\{f:$ for some ideal $I$ of $A, f $ is
an onto endomorphism of $I\}$.

 $End_5(A)=\{f:$ for some dense ideals $I,J,
\,f$ is  an homorphism from $I$ onto $J\}$.

 $End_6(A)=\langle f:$ for some ideals $I,J$ of
$A, f$ is an endomorphism from $I$ onto $J\}$.

{\bf Note} $I=A$ is allowed. All kinds of endomorphism, commute with $\cap, \cup$, preserve $0$ but not necessarily $-$. 

\item"(2)"  For $l=1\cdots, 6$ we let $Aut_l(A)=\{f \in
End_l(A):\, f(x)=0\Leftrightarrow  x=0$ for $x\in
Dom\,f\}$.
\item"(3)" For function $f,g$ whose domain is $\subseteq A$
let: $f\sim g$ if $Kerf=Kerg$ and $\{x$ 
$:f(x)=g(x)$ or both are defined not $\}$
include a dense ideal of $B/Kerf$. 

\item"(4)" Let $end_e(A)=|End_e(A)|$,\,\,
$aut_e(A)=|Aut_e(A)|$. 

Let $end^\sim_e(A)=|\{f/ \sim \, :\,f\in
End_e(A)\}|$ and $aut^\sim_e(A)=|\{f/\sim:f\in Aut_e(A)\}|$.
\item"(5)"  We allow to replace $A$ by an ideal $I$ with
 the natural changes.
\item"(6)"  We define $ Endv \sb e (A), endv_e(A)$ similarly
replacing ``onto" by ``everywhere onto" and define
$Endl \sb e(A), endl_e(A)$ similarly omitting ``onto". We
define $Endu_e(A)$  as the set of $f \in End_e(A)$ such that for every $x \in Dom\,f, f(x) \ne 0$ and ideal $I$ of $A|x$ which is dense, we have $f(x) = sup_A\{f(y): y \in I\}$. We defined naturally
$Autv_e(A), Aut l_e(A), Autv_e(A)$ etc.

{\bf Note}:  In $Endv\sb1(A)$ we mean:  for some dense
ideal $I$ of $A, y \leq f(x)\in I \Rightarrow (\exists z\leq x)
(f(z)=y)$ and  $Endl\sb1(A)=End(A)$.
 \endroster 

\mpr{\bf 3.7 Claim}. \; {\sl
\roster 
\item  $End_6(A) \supseteq End_5(A)\cup
End_4(A)\supseteq End_5(A)\cap End_4(A)\supseteq
End_3(A)\supseteq  End_2(A)$ and $End\sb2(A) \subseteq
End\sb1(A) \subseteq End\sb0(A)=EndA$.
 \item $
end_6(A)\geq max \{ end_5(A), end_4(A)\} \geq min
\{end_5(A), end_4(A)\}\geq end_3(A)\geq \mathbreak end_2(A)$, and $end\sb2(A) \leq end\sb1(A)\leq
end(A)$.
 \item  In (1) we can replace $End$ by $Aut$ or $Endv$ or
a$Endl$ or $Autv$ or $Autl$, or $Endu$ or $Autu$ and in (2) $end$ by $aut$ etc. 
We can in (2) replace $end_l$ (or $aut_l$) by $end^\sim_l$
(or $aut_l^\sim$) etc.
 \item 
$Aut_l(A)\subseteq End_l(A)$ hence $aut_l(A) \leq
end_l(A)$; and $Endv\sb e (A) \subseteq End\sb e(A) \subseteq Endl\sb e (A)$ hence $endv\sb e (A) \leq end\sb
e(A) \leq endl\sb e(A); Endv\sb e(A)=End\sb e(A)$ if
$e=0$. Also $Aut(A)=  Aut^\sim_2(A)$ hence $aut(A)= aut^\sim_2(A)$.
  \item  If $f,g,\in Endv_6(A),
\, f\sim g$ {\it then} $f,g,$ are compatible
functions,i.e. $x\in Dom f \cap Dom g \Rightarrow
f(x)=g(x)$.
\item $sub(A)\geq aut_6(A)$ and $aut(A) \le aut^\sim_3(A)$.
\item $Aut_e(A)=Autv\sb e(A)$.
\item $end\sb e(A) \leq id(A) + end^\sim\sb e(A)$ etc.
  \endroster
}

\mpi E.g.
\roster
\item"(5)" Supposse $x\in
Domf\cap Domg,\,f(x)\not=g(x)$, so $w log \, f(x)\not\leq
g(x)$ so let $z=:\,f(x)-g(x)\not=0$.  As $z\leq f(x)$ for
some $t\in Dom f$ we have $0<f(t)\leq z$, and $wlog\,
f(t)=g(t)$, and we get contradiction.
\item"(6)"  As in proof of 3.2 or see proof of 3.12
(noting $|\{I:I$ a dense ideal of $A\}| \leq sub(A)$ by 3.4A).
\endroster

\mpr{\bf 3.8 Definition}. \; For a Boolean Algebra $A$
\roster
\item  $Idc(A)=\{I:I$ an ideal of $A$, and $I^c=I\}$ where
$I^c=\{x\in A: I$ is dense below $x\}$.
\item  $idc(A)=|Idc(A)|$.
\item  $Did(A)=\{I:\, I$ a dense ideal of $A\}$.
\item  For ideals $I, J\,\, I+J$ is the minimal ideal $I$
of $A$ which include $I\cup J$ i.e. $\{x \cup y: x \in I$
and $y\in J\}$.  Similarly $\sum_{\zeta<{\xi_0}} I_\zeta$.

\item If we replace $A$ by an ideal $I$ (in 3.8
(1),(2),(3), 3.1 (2)) means we restrict ourselves to
subideals of it.
 \endroster

\mpr{\bf 3.9 Claim}. \; {\sl  For a Boolean algebra A:
\roster
 \item  $ id(A) \geq idc(A)=|comp(A)|$.
\item  $|A|\leq idc(A)=idc(A)^{\aleph_0}$ (when $A$ is
infinite, of course).
\item  If $f\in Endu_5(A)$, {\it then} $f$ has a unique
extension to an endomorphism $f^+$ of $comp(A)$
where $f^+(x)= sup_A \{f(y):y \leq x,y \in
Dom f\}$.  If $f$ is everywhere onto it is the unique
extension of $f$ in $End(comp A)$.
 \item  For $g\in End(A), \,
(\exists f \in End_5(A))f^+ \supseteq g$ {\it iff} $g \in Endu(A)$.
\item  For $f,g \in End_5(A)$ we have
$f\sim g \Leftrightarrow f^+=g^+$.
\item  For $f\in End\,u_6(A)$,
letting $a=sup_{comp(A)}\{x:\,B\restriction x \subseteq
Dom\,f\}\in comp(A)$ and $b=sup_{comp(A)}\{x: \,A\restriction x
\subseteq Rang f\} \in comp(A)$, we define $f^+\in 
HOM(comp ( A\restriction a)$, $comp (A\restriction b))$
extending $f$ by $f(x)=sup\{f(y):y\leq x, y \in Dom f\}$,
also $f^+$ is onto.
 \endroster
}
 \mpr{\bf 3.10 Claim}. {\sl \roster
\item  $id(A)\leq endl_3(A)$.
\item  $idc(A)\leq endl_3^\sim(A)$.
\item   $endl_e(A)=id(A)+endl^\sim_e(A)$ for $e=3,4,5,6$.
\item  If $f\in End_5(A)$ {\it then} $f^+ \in End_1(comp
A)$. 
\item $aut_e(A) \le id(A) + aut^\sim_m(A)$ when $e$, $m \in \{3,5\}$
or $e$, $m \in \{4,6\}$. 
 \endroster
}

\mpi
\roster
\item  For $I\in Id(A)$ let
$J_I=\{x \in B:(B\restriction x) \cap I=\{0\}\}$, so $J_I$
is an ideal, $J_I\cap I=\{0\}, \, J_I\cup I$ is dense.  Let
$F_I$ be the following map:  $Domf_I=I+J_I,\,
f_I\restriction I=id_I,\, f_I\restriction J_I=0_{J_I}$. 
Now $I\mapsto f_I$ is a one to one mapping from $Id(A)$ to
$Endl_3(A)$. 
\item  The mapping above works.
\item   Note that $id(A)\leq end_3(A)$ by part(1),
and $endl_3^\sim(A)\leq endl_3(A)$ trivially. $A$ is
infinite hence all those cardinals are infinite so
$\chi=:id(A)+endl_3^\sim(A)\leq endl_3(A)$.  The
inverse inequality is easy too.
\endroster
\par\quad (4), (5) \;  Left to the reader.

\mpr{\bf 3.11 Claim}. \; {\sl
\roster 
\item  For $x\in A$:

$\big(\exists f\in Aut(A)\big)[x\not=f(x)]$ {\it iff}

$\big(\exists f\in Aut_6(A)\big) \big(x\in Dom f \&
x\not=f(x)\big)$ {\it iff}  $(\exists y,z)[0<y \leq x \&
0<z\leq 1-x \& A\restriction y\cong A\restriction z]$.
\item  If $\langle I_\zeta:\zeta<\alpha\rangle$ is a
sequence of ideals of $A$
,$[\zeta\not=\xi\Rightarrow I_\zeta\cap I_\xi=\{0\}]$ and

$I=\sum_{i<\alpha}I_i$ {\it then}:

$id(I)\geq \pi_{\zeta<\alpha}\, id(I_\zeta)$.

$idc(I)\geq\pi_{\zeta<\alpha}\, idc(I_\zeta)$.

$aut(I)\geq\pi_{\zeta<\alpha}\, aut(I_\zeta)$.

$end(I)\geq \pi_{\zeta<\alpha}\, end(J_\zeta)$.

Similarly for $end_l, aut_l, end^\sim_l,
aut^\sim_l$ etc. 
\endroster
There are many more easy relations, but for our aim the
main point is }

\mpr{\bf 3.12 Main Lemma}. \; {\sl  For an infinite Boolean
Algebra A:
 \roster
\item  $aut^\sim_e(A)$ for $e=3,5$ are equal or both
finite (and we can restrict ourselves to automorphisms of
order 2).
\item  If  for some $e \in \{3,4,5,6\}$ we have
$aut^\sim_e(A)> idc(A) $ {\it then}
$aut^\sim_e(A)$ for $e=3,4,5,6$ are all equall.
 \item 
$aut_3^\sim(A)=aut^\sim_3(A)^{\aleph_0}$ or $aut_3(A)$ is
finite. 
\item $autv^\sim\sb e(A)=autv\sb3(A) +idc(A)$
for $e=4,6$.
\endroster 
}

\mpi  Let $J=:\{x\in A:$ for every $f\in
Aut(A),\,f\restriction  (A\restriction
x)=id_{A\restriction x}\}$.

The function $F_{1,\cdots}, F_5$ satisfying $y\leq
x\Rightarrow F_e(y)\leq F_e(x)$ are functions from $A$ to
$ord$ defined below; and we let: 

$K=:\{x\in A:$

\qquad\qquad (i) for some $y, x\cap y=0$ and
$A\restriction x\cong A\restriction y$.

\qquad\qquad (ii) for $e=1,\cdots, 5$ and $0<y<x$ we have
$F_e(y)=F_e(y)\}$.

Where:

$F_1(x)=$ the cardinality of $A\restriction x$.

$F_2(x)=idc(A\restriction x)$.

$F_3(x)=aut(A\restriction x)$.

$F_4(x)= aut_3(A\restriction x)$.

$F\sb5(x)=aut^\sim_3(A|x)$.

Now

$(*)_1 \,\, J$ is an ideal of $B$.

$(*)_2 \,\,K$ is downward closed.

$(*)_3 \,\, K\cup J$ is dense.

Choose $\{x_i:i<\alpha\}$ maximal such that:
\roster
\item"(a)"  $x_i\in K, \,x\sb i>0$,
\item"(b)"  if $i\not=j,\, 0<y'\leq x_i,\, 0<y''\leq x_j$
{\it then} $A\restriction y' \not\cong A\restriction y''$.
\endroster

Let $K_i=\{y:$ for some $y'\leq x_i, \, A\restriction y
\cong A\restriction y'\}$, and $K_i^+:$ the ideal $K_i$
generate.  Now

$(*)_4\,\,  \bigcup_{i<\alpha} K_i$ is dense in $K$
[hence $J\cup\bigcup_{i<\alpha}K_i$ is dense in $A$].

$(*)_5$\ \ For $i\not=j$ we have $K_i\cap K_j=\{0\},
\,K_i\cap J=\{0\}$.

Clearly,

$(*)_6$ \,\, for $f \in Aut_6(A)$ we have:

\qquad\qquad (a) \; $f$ is the identity on $J\cap Dom
f$. \par
\qquad\qquad (b)\;\;  for $x\in Domf$ we have  $x\in
K_i\Leftrightarrow f(x)\in K_i$.

So

$(*)_7\,\,
aut^\sim_e(A)=\pi_{i<\alpha}\,aut^\sim_e(K^+_i)$ for
$e=3,5$ and $aut^\sim_e(A)=idc(J) \times \pi_{i<\alpha}
aut^\sim_e(K^+_i)$ for $e=4,6$

We shall prove:

$(*)_8$  For each $i$, \ one of the following ocurrs:  
\roster 

\qquad\qquad \item"(a)"  $aut^\sim_e(K^+_i)$ for
$e=3,4,5,6$ are all finite $>1$,\par
\qquad\qquad \item"(b)" For some infinite $\kappa$ we have
$aut^\sim_e(K^+_i)= autv^\sim\sb
e(K^+\sb i)^\kappa \geq idc(K_i)$  for $e=3,4,5,6$

(really we can use $F\sb6(i)= sup\{\kappa^+: B
\restriction x$ has $\kappa$ pairwise disjoint non zero
members$\}$ and any such $\kappa$ is OK for (b)).
 \endroster
 
\mpr{\bf Case 1}.  $x_i$ is an atom.

This is clear:  let $\lambda_i=:|K_i|$, if it is
infinite, $aut^\sim_e (K^+_i)=2^{\lambda_i}$ for
$e=3,4,5,6$ so case (b) in $(*)_8$ ocurrs.

If $\lambda_i$ is finite, $1<aut^\sim_e (K^+_i)<\aleph_0$
(we can compute exactly), so case (a)   in $(*)_8$.  In
fact in all cases we can use just automorphism of order 2.

So we can assume

\mpr{\bf Case 2 }. \;  not Case 1, so $B
\restriction x_i$ is atomless, hence $idc(B\restriction
x_i) \geq 2^{\aleph_0}$

Let $\langle J_{i,\zeta}, J'_{i,\zeta}:\zeta
<\zeta_i\rangle$ be a sequence such that:
\roster
\item"$(\alpha)$"  $J_{i,\zeta}, \,J'_{i,\zeta}$
are ideals $\subseteq K^+_i$ and $\not=\{0\}$,
\item"$(\beta)$"  $\{0\}\not= J_{i,\zeta} \subseteq
A\restriction x_i$,
\item"$(\gamma)$"  $J_{i,\zeta}\cong J'_{i,\zeta}$, an $h_{i,\zeta}$ an isomorphism from $J_{i,\zeta}$ onto $J'_{i,\zeta}$,
\item"$(\delta)$" $\wedge_{\zeta<\xi<{\zeta_i}} \,
J'_{i,\zeta} \cap J'_{i,\xi}=\{0\}$,
\item"$(\epsilon)$" if $y \in K^+_\zeta$ is
disjoint to all members $\cup\sb{\xi\le\zeta}J'_{i,\xi}$
of then for some $y' \le y$ and $z \in J_{i,\zeta}, B|z' \cong B|y'$ hence
\item"$(\zeta)$" if $\zeta <\xi <\zeta\sb i$ then
$J\sb{i,\xi} \cap J\sb{i,\xi}$ is a dense subset of
$J\sb{i,\xi}$ 
 \endroster
Now

$(*)_9\,\, \bigcup_{\zeta<\zeta_i}J'_{i,\zeta}$ is a dense
subset of $K_i^+$.

$(*)_{10}\,\, \zeta_i\geq 2$ (as there is $y,x_i\cap y=0$
and $A\restriction y\cong A\restriction x_i$ as $x_i\in
K)$.

$(*)_{11} \,\,
idc(K^+_i)=\pi_{\zeta<\zeta_i}idc(J'_{i,\zeta})$.

By the definition of $K$ (see choice of $F_2)$, we have
$0<y\leq x_i\Rightarrow idc(A\restriction y)=(A\restriction
x)$

hence $idc(J'_{i,\zeta})=idc(A\restriction x_i)$ so

$(*)_{12} \,\, idc(K^+_i)=[idc(A\restriction
x_i)]^{|\zeta_i|}$ hence $[idc
(K^+_i)]^{|\zeta_i|}=idc(K^+_i)$.

Easily 

$(*)_{13}$ $aut_6^\sim( K^+_i)\geq \,
 aut_3^\sim(K^+_i)\geq idc(A\restriction
x_i)$.

(and even by automorphism of order 2).

 Last inequality:
for each $z\in Idc(A\restriction x_i)$ there is
$z', z'\cap x_i=0$, and $A\restriction z\cong
A\restriction z'$ and let $g$ be such isomorphism, let $I_z$ be the
ideal of $A$ generated by $\{x \in K^+_i: x \le z$ or $x \le z'$ or $x
\cap z = x \cap z' = 0\}$  
 $g_z\in Aut(I_z) \subseteq Aut_3(K^+_3), g_z(y)=:(y-z\cup z')\cup
g(y\cap z)\cup g^{-1}(y\cap z')$.

\mpr{\bf Case  2A}. 
$\zeta_i<\aleph_0$ (and not Case 1)

Let for $n < \omega$ $x\sb{i,n}\leq x\sb i, x\sb{i,n}\not=
0$, and $[n\not=m\Rightarrow  x\sb{i,n} \cap x\sb{i,m}=0]$, and
$I\sb{i,\omega}=\{y\in K\sb i^+: \wedge\sb n y \cap
x\sb{i,n}=0 \}$ and $I\sb{i,n}=\{z:z\leq x\sb{i,n}\}:$ So
$aut_3^\sim(K^+\sb i) \geq aut_3^\sim (\sum_{\alpha \le \omega} I\sb{i,\alpha})
\geq \pi_{n < \omega} aut_3^\sim (B\restriction
x\sb{i,n})=aut_3^\sim(B\restriction x\sb i)^{\aleph \sb0}$. 

Similarly to the proof of Case 2B below (but easier) we can show that Case (b) of
$(*)_8$ ocurrs.  From now on we assume 

\mpr{\bf  Case 2B}. \;  Not Case 2A, so
$\zeta_i\geq\aleph_0$ 
 
For every $f\in Aut_6(K^+_i)$ let, for
$\zeta,\xi<\zeta_i$:

$L^f_{\zeta,\xi}=\{x\in J'_{i,\zeta}:\, f(x)\in
J'_{i,\xi}\}$.

$M^f_{\zeta,\xi}=\{f(x):\, x\in J'_{i,\zeta},\, f(x)\in
J'_{i,\xi}\}$.

So the number of possible $ \bar L^f =\langle
sup L^f_{\zeta,\xi}:\, \zeta\not=\xi < \zeta_i
\rangle$ and $\bar M^f=\langle sup (M^f_{\zeta,\xi}):\,
\zeta\not=\xi<\zeta\rangle \,\text{is} \leq 
[idc(A\restriction x_i)]^{|\zeta_i|}$ and for fixed
$\bar L, \bar M$ the number of $f\in
aut_6^\sim(K^+_i)$ for which $\bar L^f=\bar L,
\bar M^f=\bar M$ is $\leq \pi_{\zeta,\xi} |\{f/\sim$ 
$: f$ an isomorphism from a dense subset of $L_{\zeta,\xi}$
onto a dense subset of $M_{\zeta,\xi} \}| \leq
\pi_{\zeta,\xi}aut^{\sim}_3(L_{\zeta,\xi}) \leq
\pi_{\zeta,\xi} \,aut_3^\sim(A\restriction x_i)=
[aut_3^\sim (A\restriction x_i)]^{|\zeta_i|})$.

(In the last equality we use $F_5$ in the definition of $K$: for the
last $\leq$, note we can replace $L_{\zeta,\xi}$ by isomosrphic ideal
$L^*_{\zeta,\xi}$ of $A\restriction x_i$, and letting
$L^-_{\zeta,\xi}=\{y\;:\: y\in A\restriction x_i, \;(\forall z\in
L^*_{\zeta,\xi})[y\cap z=0]\}$ we can extend every $f\in
Aut_3(L^*(\zeta,\xi)$ to $f'\in Aut_3(A\restriction x)$ by letting
$f' (x\cup y)=f(x)\cup y$ for $x\in Dom(f)$, $y\in L^-_{\zeta,\xi}$)

So 

$(*)_{14} \ \ aut_6^\sim(K^+_i) \leq
[aut_3^\sim(A\restriction x_i) + idc(A\restriction
x_i)]^{|\zeta_i|}$.

Now for each $\xi$ such that $2\xi+1<\zeta_i$ we can
choose  $y_{i,2\xi} \in J_{i, 2\xi}, y_{i,2\xi+1} \in
J_{i,2\xi+1}$ such that $A\restriction
y_{i,2\xi}\cong A\restriction y_{i,2\xi+1}$ and
$\langle g_{i,\xi,\alpha}: \alpha
<aut_3^\sim(A\restriction x_i) \rangle$ such that
$g_{i,\xi,\alpha}$ is an isomorphism from a dense subset of
$A\restriction y_{i,\xi,\alpha}$ onto a dense subset of
$A\restriction y_{i,2\xi+1}, \langle 
g_{i,\xi,\alpha}^+:\alpha < aut^{\sim}_3(A\restriction x_i) \rangle$
pairwise distinct. Let $\zeta^*_i$ be minimal such that $2\zeta^*_i \le \zeta_i$.

Now for every sequence
$\bar \alpha=\langle\alpha_\xi:\xi<\zeta^*_i\rangle, \,
\alpha_\xi < aut_3^\sim(A\restriction x_i)$, we define
$g_{i,\bar \alpha} \in Aut_3(K^+_i)$ of order two (see condition $(\gamma)$):
$$g_{i,\bar\alpha} \restriction (y_{i,2\xi}
\cup y_{i,2\xi+1})=h_{2\xi+1}\circ g_{i,\xi,\alpha} \circ
h^{-1}_{2\xi} \cup h_{2\xi} \circ g^{-1}_{i\xi,\alpha}\circ h^{-1}_{2\xi+1}$$ 

and if $y\in K^+_i \backslash \{0\}$ and
$\wedge_{\xi<\zeta^*_i} (y\cap y_{i,\xi}=0 \rangle $
 {\it then} $g_{i,\bar\alpha}(y)=0$.

Lastly if $y$ is the disjoint union of $y_0, \cdots, y_n$ and each
$g_{\bar \alpha}(y_e)$ was defined then we define $g_{\bar \alpha}(y_0
\cup \cdots \cup y_n) = g_{\bar \alpha}(y_0) \cup \cdots \cup g_{\bar
\alpha}(y_n)$. The reader may check that $g_{\bar \alpha} \in Aut_3(K^+_i)$.   

The mapping $\bar\alpha\mapsto
g_{i,\bar\alpha}$ show:

$(*)_{15}\, \, aut_3^\sim(K^+_i) \geq
aut_3^\sim(A\restriction x_i)^{|\zeta_i|}$. 

Easily, choosing $y_i\in J_{i,1}$ (possible as $\zeta_i\geq
2$)

$(*)_{16}\, \, aut_3^\sim(K^+_i)\geq aut_3^\sim(A|(x_i\cup y_i))\geq
idc(A\restriction y) = idc(A\restriction x_i)\geq \aleph_0$.

 Together  by $(*)_{13}, (*)_{14},
(*)_{15},(*)_{16}$ for $e=3, 6$:
$$aut_e^\sim(K^+_i)=aut_3^\sim(A\restriction
x_i)^{|\zeta_i|}.$$

By 3.7(1) this holds for $e=4, 5$, hence $(*)_8$ has been 
 proved.

Now by $(*)_7 + (*)_8$ the four parts of 3.12 follows.

\mpr{\bf 3.13 Conclusion} :  $aut(A)^{\aleph_0}
\leq \, sub(A)$, also $aut\sb e(A)^{\aleph\sb 0}\leq
sub(A)$ for $e=3,4,5,6$.

Note:  even if $aut(A)$ is finite, $A$ infinite, still
$2^{\aleph_0}\leq \,sub(A)$ (for $A$ infinite).

\mpi  By 3.7(3) + 3.7(6) $aut(A)=aut^\sim_0(A)
\leq id(A) + aut^\sim_3(A) \leq  id(A) + Pend(A)
\leq sub (A)$ but by S. Shelah [6] $id(A)^{\aleph\sb0}=id(A) $
and  by 3.12(3) $aut^\sim\sb3
(A)^{\aleph\sb0}=aut^\sim\sb 3(A)$, together we can
finish the first inequality, the second is similar using 3.10(5).

\mpr{\bf 3.14 Claim} : \; {\sl For an (inifinite) Boolean
algebra $A$ we have:
\roster
\item    
$end_e(A)=id(A)+aut_3(A)$ for $e=3,4,5,6$.
\item $end_6(A)^{\aleph_0}\leq\, sub(A)$.
\item $end_e(A)^{\aleph_0}\leq\,sub(A)$ for
$e=2,3,4,5,6$.
 \endroster
}

\mpi
\roster
\item Clearly $end_3(A) \geq\,aut_3(a)$ and
$end_3(A)\geq\,id(A)$ (as the mapping in the proof of
 3.10(1) examplify).

On the other hand we can attach to every $f \in
End_6(A)$ three
ideals $I_1(f)=Dom\,f ,\, I_2(f)=Rang(f)$
and $I_3(f)=Ker(f)$.  Now the number of triples
$\bar I$ of ideal of $A$ has cardinality $id(A)$ and
for each such $\bar I$:

$ \{f\in End_6(A):\, I_e(f)=I_e$ for
$e=1,2,3\}$
has cardinality $|Aut(I_2)|$ which is $\leq \,aut_3(A)$. 
By 3.7(2) we can finish.
 \item  Remember also
$id(A)^{\aleph_0}=id(A)$ by [6] and part (1) and 3.13.
\item  By part (2) and 3.7. 
\endroster  

 \mpr{\bf  \S 4 The width of the Boolean algebra
}. 
 
\mpr{\bf 4.1 Definition}. \;   For a Boolean
algebra $B$ let: 
(1) $A \subseteq B$ is an antichain if $x \in A \& y \in
A\& x \not= y \Rightarrow x \not\leq y$ (i.e. $A$ is a
set of pairwise incomparable elements).

(2) Width of $B$, $w(B)$ is $sup\{|A|: A\subseteq B$
antichain$\},  w^+(B)=\cup\{|A|^+:A
\subseteq B$ antichain$\}$

{\rm E. C. Milner and M. Poizat [3], answering a question of E. K. van
Dowen, D. Monk and M. Rubin [7]}
{\rm proved $cf(w^+(B))\not=
\aleph_0$.}

{\rm In S. Shelah [5]  we claim: if
$\lambda>cf\lambda>\aleph_0$, for some generic extension
of the universe preserving cardinalities and cofinalities,
for some $B, w^+(B)=\lambda$.  We retract this and replace
it by the theorem 4.2 below.}

{\rm For weakly inaccessibles we still have the
consistency.  Moreover, if $\lambda$ is a limit
uncountable regular cardinal, $S\subseteq \lambda$
stationary not reflecting and $\diamond_S$ {it then} we
have such an example for $\lambda$.}

\mpr{\bf 4.2 Theorem}. \; {\sl  For an infinite Boolean algebra
$B, w^+(B)$ is an uncountable regular cardinal.
}

\mpi  As $B$ is infinite it has an antichain
$A, |A|=\aleph_0$, [if $B$ has finitely many atoms clear,
if not it has a subalgebra which is atomless, without loss
of generality countable and check].  So
$\lambda=:w^+(B)>\aleph_0$.  Assume $\kappa=: cf\lambda
<\lambda$; let $\lambda=\sum_{i<\kappa} \lambda_i,\,
cf\lambda +\sum_{j<i}\lambda_j < \lambda_i<\lambda$ and
let $A_i\subseteq B$ be an antichain of cardinality
$\lambda_i^+$ ( exist by the choice of $\lambda$).  Let
$A=\bigcup_{i<\kappa} A_i$, so $|A|=\lambda$.  Choose such
$\langle A_i:\,i<\kappa\rangle$ such that, if possible

(*)$i<j<\kappa, \, x\in A_i,\,
y\in A_j\Rightarrow y\not\leq x$.

For $x\in B$ let $A[>,x] =\{y \in A:\,y>x\},  A[<,x]=\{y
\in A:y<x\}$,

 $A[>,\mu ]=\{x\in
A:|A[>,x]|<\mu\} , A[<,\mu\rbrack =\{x\in
A:|A[<,x] |<\mu\} $.

\mpr{\bf Case 1}. \;  For some $\mu <\lambda,
\,A[>,\mu] $ has cardinality $\lambda$.

By Hajnal free subset theorem, there is a set $E
\subseteq A[>,\mu]$ of cardinality $\lambda$ such
that:

 $x \not= y \, \& \, x \in E \, \& \, y \in E \Rightarrow x
\not\in A[>,y] \,  \& \, y \not\in A[>,x ]$. So $E$ witness
$w^+(B)>\lambda$. 

\mpr{\bf Case 2}. \;  For some $\mu < \lambda, \,
A[<,\mu \rbrack $ has cardinality $\lambda$.

Same proof.

\mpr{\bf Case 3}. \;  For every $ i<\kappa$ there
is $x \in A$ such that $\lambda_i <|A[>,x]
|<\lambda$.

Let for $i<\kappa,\,x_i\in A$ be such that
$\lambda_i<|A[>,x_i]| <\lambda$.  Let $u \subseteq
\kappa$ be such that: $\kappa= sup u$ and for $i \in
u, \, \lambda_i> \sum_{j\in u \cap i}|A [>,x_i]|$ (choose
the members of $u$ inductively).  By renaming without
loss of generality $u=\kappa$. 
Clearly $A[>,x_i\rbrack \backslash\cup_{j<i}
A[>,x_j\rbrack$ has cardinality $ > \lambda_i$.

As $\lambda_i>\kappa$ (by its choice) and
$A=\cup_{j<\kappa} A_j$, clearly for each $i$ there is
$\alpha(i) <\kappa$ such that $( A[>,x_i]\backslash
\cup_{j<i} A[>,x_j]) \cap
A_{\alpha(i)}$ has cardinality $>\lambda_i$;
necessarily $\alpha(i)\geq i$.

For some unbounded $u \subseteq \kappa$ we have $[i\in u \&
j \in u \&i< j\Rightarrow \alpha(i)<j]$; without loss of generality
$u=\kappa$, $\alpha_i=i$. Let
$A^*_i$ be a subset of $( A[>,x_i]\backslash
\cup_{j<i} A[<,x_j]) \cap
A_{\alpha(i)}$ of cardinality $\lambda_i^+$.  Now
$\langle A_i^*:i<\kappa\rangle$ satisfies: $A^*_i\subseteq
B$ is an antichain of cardinality $\lambda^+_i$ and 
\roster 
\item"$(*)'$" $i<j,\, x \in A_i^*,  y\in
A^*_j\Rightarrow x\not\leq y$ \qquad (otherwise
$x_i\leq x\leq y\not\in A[>,x_i]$, contradiction).
\endroster
So $\langle A'_i=:\{1_B-x:x
\in A^*_i\}:i<\kappa\rangle$ satisfies $A'_i\subseteq
B$ is an antichain of $B$ of cardinality  $\lambda^+_i$
and also (*) above (check).  So by the choice of $\langle
A_i:i<\kappa\rangle$, it satisfies (*).  By $(*)+(*)'$,
$A^*=\cup_{i<\kappa} A^*_i$ is an antichain of $B$ of
cardinality $\lambda$, so $w^+(B)>\lambda$.

\mpr{\bf Case 4}. \; For every $i < x$ there is $x \in A$ such that
$\lambda_i < |A[<, x]| < \lambda$. 

Similar to Case 3.

\mpr{\bf Case 5}. \;   None of the previous cases.

By ``not Case 3" for some $i(*)<\kappa$, for no $x\in A$
is $\lambda_{i(*)} <| A  [>,x]| <\lambda$.  By
not Case 2,  $A[<,\lambda^+_{i(*)}]$ has
cardinality $<\lambda$.  By not Case 1
$A[>,\lambda^+_{i(*)}]$ has cardinality
$<\lambda$.  

Choose $x^*\in  A\backslash
A[<,\lambda^+_{i(*)}]\backslash A[>,\lambda^+_{i(*)}]$ so
$A[>,x^*]$  has cardinality $\geq
\lambda_{i(*)}^+>\lambda_{i(*)}$, hence by
the choice of $i(*)$ we have  
$A[>,x^*]$ has cardinality $\lambda$.

As $\lambda_{i(*)} >\kappa$, for some $j(*),
A[<,x^*] \cap A_{j(*)}$ has cardinality
$>\lambda_{i(*)}$, so choose distinct $y_i\in
A[<,x^*] \cap A_{j(*)}$ for $i<\kappa$.  
Now $y_i<x^*$ (as $y_i\in A[<,x^*]$), and
$[i\not=j\Rightarrow y_i\not\leq y_j]$ (as they are
distinct and in $A_{j(*)}$).

Let $A'_i=A_i \cap A[>,x^*]$, so $A'_i\subseteq
A_i$ hence is an antichain of $B$, and 
$$\cup_i A'_i=(\cup A_i) \cap
A[>,x^*]=A \cap A[>,x^*]=A[>, x^*]$$
So each $A'_i$ is an antichain, its member are $> x^*$ and
$|\cup_i A'_i|$ is $\lambda$ as $|A[>,x^*\rbrack|$ is.

Now
$$A'=\cup_{i<\kappa} \{y_i\cup(x-x^*):x\in
A'_i\}$$
is an antichain of $B$ of cardinality $\lambda$, so
$w^+(B)>\lambda$, as required.

 \enddocument